\documentclass[a4paper, 
11pt
]{article}

\newcommand{\beginpipe}[2]
{
  \draw (#1    , 0  ) --
        (#2-0.5, 0  ) --
        (#2,     0.8) --
        (#2-0.5, 1.6  ) --
        (#1,     1.6  ) --
        (#1,     0  )

}

\newcommand{\pipe}[3]
{
    \filldraw[color=black, fill=#3] (#1-0.5, 0) --
        (#2-0.5, 0) --
        (#2,     .8) --
        (#2-0.5, 1.6) --
        (#1-0.5, 1.6) --
        (#1,     .8) --
        (#1-0.5, 0)

}

\newcommand{\epipe}[2]
{
  \draw (#1-0.5, 0  ) --
        (#2,     0  ) --
        (#2,     1.6  ) --
        (#1-0.5, 1.6  ) --
        (#1,     0.8) --
        (#1-0.5, 0  )

}

\usepackage{subfig}
\usepackage{pgfplots}
\usepackage{pgfplotstable}
\usepackage{mathtools}
\usepackage{multicol}
\usepackage{comment}
\usepackage{booktabs}
\pgfplotsset{compat=1.5}
\usepackage{amssymb}
\usepackage{url}

\usepackage{pdfsync}
\usepackage{float}
\usepackage{tabularx}
\usepackage{enumerate}
\usepackage{array}
\usepackage{xspace}

\usepackage{authblk}

\renewcommand{\u}{\ensuremath{\boldsymbol{u}}}

\newcommand{\p}{\ensuremath{p}}

\newcommand{\n}{\ensuremath{\boldsymbol{n}}}
\newcommand{\x}{\ensuremath{\boldsymbol{x}}}
\newcommand{\mupar}{\ensuremath{\boldsymbol{\mu}}}
\newcommand{\etapar}{\ensuremath{\boldsymbol{\eta}}}

\begin{document}

\title{Combined parameter and model reduction of cardiovascular problems by means of active subspaces and POD-Galerkin methods}

\author[]{Marco~Tezzele\footnote{marco.tezzele@sissa.it}}
\author[]{Francesco~Ballarin\footnote{francesco.ballarin@sissa.it}}
\author[]{Gianluigi~Rozza\footnote{gianluigi.rozza@sissa.it}}

\affil[]{Mathematics Area, mathLab, SISSA, via Bonomea 265, I-34136 Trieste, Italy}

\maketitle

\begin{abstract}
In this chapter we introduce a combined parameter and model
  reduction methodology and present its application to the efficient
  numerical estimation of a pressure drop in a set of deformed
  carotids. The aim is to simulate a wide range of possible occlusions
  after the bifurcation of the carotid.
A parametric description of the admissible deformations, based on radial basis functions interpolation, is introduced.
Since the parameter space may be very large, the first step in the
combined reduction technique is to look for active subspaces in order to reduce the parameter space dimension. 
Then, we rely on model order reduction methods over the lower
dimensional parameter subspace, based on a POD-Galerkin approach, to
further reduce the required computational effort and enhance
computational efficiency.
\end{abstract}

\section{Introduction}
\label{sec:intro}
Numerical simulations of biomedical problems is a topic of large interest nowadays, especially for what concerns the application of shape optimization techniques aimed at the improvement of long-term outcomes of clinical interventions \cite{Marsden2014,MSA5}. Several challenging aspects can be identified in such a task, especially when seeking a personalized (patient-specific) treatment~\cite{Wang2012,gonzalez2016computational,brown2015building}: model construction and segmentation, numerical solution of the underlying fluid dynamics equations, assimilation of clinical data (e.g. for boundary conditions), choice of the cost functional and medical indices to be optimized~\cite{AgoshkovQuarteroniRozza2007,AgoshkovQuarteroniRozza2006}, as well as model deformation during the optimization procedure. The latter is a topic of remarkable interest, since it is well known that local geometrical features may severely affect the computational fluid dynamics (CFD) simulation and thus the results of the optimization \cite{DoorlySherwin2009}.

A challenge in the applicability of optimization procedures in the clinical environment is the large computational time that each procedure would entail. Indeed, such optimal control problems are usually tackled by means of iterative solvers that require several expensive CFD simulations for different geometrical configurations \cite{gunzburger2003perspectives}. To this end, several authors have proposed to employ computational reduction techniques based on reduced order (or surrogate) models. We refer to \cite{manzoni2012model,McLeodCaiazzoFernandezMansiVignon-ClementelSermesantPennecBoudjemlineGerbeau2010,guibert2014group,CuetoChinesta2014,BallarinJCP,Ballarin2017} for a few representative applications, as well as e.g.\ to \cite{hesthaven2015certified,chinesta2013proper} for introductory textbooks on the underlying methodology. These methods rely on the definition of a parameter space, which is related to the set of admissible deformations that can be considered during the iterative optimization procedure, and its exploration to retain the most relevant features of the CFD solution on the parameter space.
It is however well known that reduced order methods suffer from the so-called ``curse of dimensionality'' if the parameter space is high-dimensional. Although there exists techniques to account for high-dimensional parameter spaces in the reduced order modelling framework, based e.g.\ on sparse grids \cite{Chen2014SparseGrid,luca} or on a proper weighting of the parameter space \cite{Chen2014Weighted,davide}, the approach that we propose in this manuscript is different and aims at reducing the high-dimensional parameter space as well, while preserving a very broad set of admissible deformations.

We have recently dealt with several techniques concerning efficient shape parametrization techniques in the framework of reduced order modelling. A first possible choice is related to the shape morphing method itself. For our goals it suffices to classify them in two groups: \emph{general purpose} or \emph{problem specific}. 
The design of a \emph{problem specific} shape parametrization technique should aim at reducing the high-dimensional parameter space; for instance, the centerlines-based approach proposed in \cite{BallarinJCP} is able to reduce the parameter space accounting only for deformations in a cylindrical coordinates frame of reference. In contrast, in this chapter we show a parameter space reduction technique for \emph{general purpose} shape morphing methods.
Among possible \emph{general purpose} methods we mention \emph{Free Form Deformation} (FFD) \cite{sederbergparry1986,LassilaRozza2010}, \emph{Radial Basis Functions} (RBF) interpolation \cite{buhmann2003radial,morris2008cfd,manzoni2012model} or \emph{Inverse Distance Weighting} (IDW) interpolation \cite{shepard1968,witteveenbijl2009,forti2014efficient}. Broadly speaking, the aforementioned methods require the displacement of some control points to induce a deformation on the domain, and we identify the parameters as the displacements of the control points.
In these context, our earlier approaches to parameter space reduction have relied on screening procedures based on Morris' randomized one-at-a-time design \cite{Morris,BallarinManzoniRozzaSalsa2013}, modal analysis \cite{forti2014efficient} or semi-automatic reduction of the number of control points \cite{DAmario}. Each of these approaches results in lower dimensional parameter spaces by retaining an ``optimal'' subset of the possible control points. Unfortunately, when requiring very low dimensional parameter spaces (order of unity) the resulting set of admissible deformations may be considerably shrunk. Not only, another goal to achieve for complex shape parametrizations is also an efficient positioning of control points, requiring versatility and ``full'' capability in the geometric representation. To this end, in this manuscript we propose to exploit an active subspaces (AS) method \cite{constantine2015active}, building on our previous experience on a shape optimization of a naval engineering problem in \cite{tezzele2017dimension}. AS has been employed in many real world problems; among others, we mention aerodynamic shape optimization~\cite{lukaczyk2014active}, integrated hydrologic model~\cite{jefferson2016reprint}, the parameter reduction
for the HyShot II scramjet model~\cite{constantine2015exploiting}, a
satellite system model~\cite{hu2016discovering}. 
The main difference between our previous approaches and the AS property is that our former approaches were constraining the search of a lower dimensional parameter subspace to be parallel to a subset of the axes of the parameter space, while the AS method will automatically identify the ``optimal'' lower dimensional subspace without any such constraint, taking a linear combination of all the original parameters. To show an example of the resulting methodology we will consider a cardiovascular test case related to the computation of the pressure drop of a series of deformed carotids.

\begin{figure}[t]
\centering
\begin{tikzpicture}[node distance = 2.5cm, auto]
	\beginpipe{0}{2.5};
	\pipe{2.55}{4.3}{white};
	\pipe{4.35}{6.0}{white};
	\pipe{6.05}{7.7}{white};
	\pipe{7.75}{9.5}{white};
        \epipe{9.55}{12.15};
        
	\node[align=left] at (1.0,.8) {Problem\\settings\\(sec. 2)};
	\node[align=left] at (3.3,.8) {RBF\\(sec. 3)};
	\node[align=left] at (5.1,.8) {FEM\\(sec. 4)};
	\node[align=left] at (6.75,.8) {   AS\\(sec. 5)};
	\node[align=center] at (8.45,.8) {POD\\(sec. 6)};
	\node[align=center] at (10.7,.8) {Reduced CFD\\evaluation\\(sec. 7)};
\end{tikzpicture}
\caption{Outline of the chapter}
\label{fig:pipeline}
\end{figure}
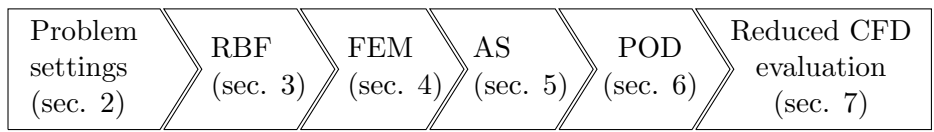

The outline of the chapter is presented in Figure~\ref{fig:pipeline}. The formulation of the problem of estimation of the pressure drop across a stenosed carotid artery is summarized in Section~\ref{sec:problem}. Shape morphing by means of RBF interpolation is then introduced in Section~\ref{sec:rbf}. The high fidelity method, based on finite elements, is briefly summarized in Section~\ref{sec:solver}. The computed values of the quantity of interest will be used to train the AS reduction of the parameter space in Section~\ref{sec:as}.
The same high fidelity solver will be used in Section~\ref{sec:pod} to train a Proper Orthogonal Decomposition (POD)-Galerkin method on the lower dimensional parameter subspace. This combination will further enhance computational efficiency for the procedure. Numerical results and error analyses of the whole pipeline will be presented in Section~\ref{sec:results}. Conclusions and future perspectives follow in Section~\ref{sec:the_end}.

\section{A model cardiovascular problem: pressure drop estimation across a stenosis}
\label{sec:problem}
In this section we introduce the problem of the estimation of the pressure drop across two parametrized stenoses in a carotid bifurcation.

\begin{figure}[t]
\centering
\includegraphics[height=0.5\textwidth]{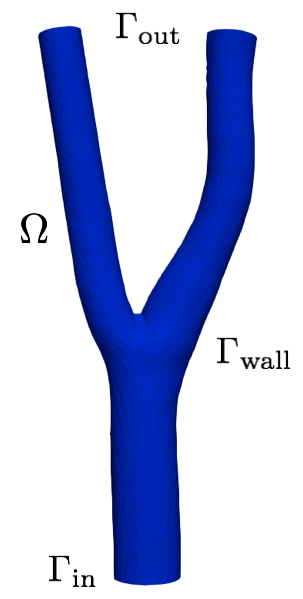}
\caption{Representation of the reference domain $\Omega$, inlet section $\Gamma_{\text{in}}$, rigid walls $\Gamma_{\text{wall}}$ and outlet section $\Gamma_{\text{out}}$.}
\label{fig:domain}
\end{figure}

Let $\Omega \subset \mathbb{R}^n$, $n = 3$, be a domain (see Figure \ref{fig:domain}), obtained from the INRIA 3D Meshes Research Database \cite{mesh}, that describes an idealized carotid bifurcation. We will call $\Omega$ the \emph{reference} domain; for practical reasons this domain happens to correspond to the healthy case (no stenoses), even though this assumption is not fundamental for the remainder of the paper.

Let $\mathbb{D} \subset \mathbb{R}^m$, be the set of parameters, that we assume to be a box in $\mathbb{R}^m$, for $m \in \mathbb{N}$.
Moreover, let $\mathcal{M}(\x; \mupar): \mathbb{R}^n \to \mathbb{R}^n$, with $\mupar~\in~\mathbb{D}$, be a shape morphing that maps the reference domain $\Omega$ into the deformed domain $\Omega(\mupar)$ as follows:
\begin{equation*}
\Omega(\mupar) = \mathcal{M}(\Omega; \mupar).
\end{equation*}
We refer to Section \ref{sec:rbf} for the actual definition of $m$ and $\mathcal{M}$ for the case at hand.

\begin{figure}[t]
\centering
\includegraphics[height=0.75\textwidth]{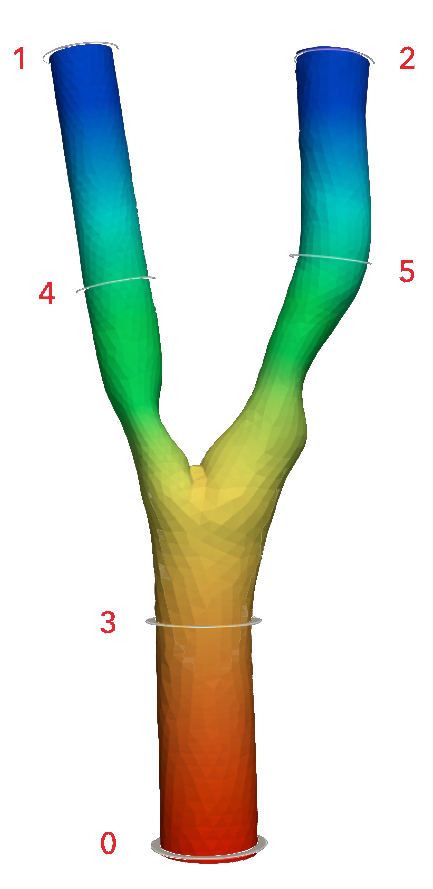}
\caption{Sections $S_0, \hdots, S_5$ employed in the relative pressure drop computation.}
\label{fig:output_sections}
\end{figure}

Let us consider the following steady Navier-Stokes equations: for any $\mupar~\in~\mathbb{D}$, find $(\u(\mupar), \p(\mupar)) \in H^1(\Omega(\mupar);\;\mathbb{R}^n) \times L^2(\Omega(\mupar))$ such that:
\begin{equation}
\begin{cases}
- \nu  \Delta \u(\mupar) + \u(\mupar)  \cdot \nabla \u(\mupar) + \nabla \p(\mupar) = \boldsymbol{0} & \text{ in } \Omega(\mupar),\\
\text{div}\ \u(\mupar)   = 0 & \text{ in } \Omega(\mupar),\\
\u(\mupar)   = \boldsymbol{u}_{\text{in}} & \text{ on } \Gamma_{\text{in}},\\
\u(\mupar)   = \boldsymbol{0}, & \text{ on } \Gamma_{\text{wall}}(\mupar),\\
\nu \displaystyle \frac{\partial \u(\mupar) }{\partial \n} - \p(\mupar) \n = \boldsymbol{0}, & \text{ on } \Gamma_{\text{out}},\\
\end{cases}
\label{eq:ns}
\end{equation}
where $\Gamma_{\text{in}}$ is the inlet section, $\Gamma_{\text{out}}$ the outlet section and $\Gamma_{\text{wall}}(\mupar)$ are (rigid) walls of the carotid artery. Since our interest is to vary the degree of stenosis immediately after the bifurcation point (see Figure \ref{fig:output_sections}), we assume that $\Gamma_{\text{in}}$ and $\Gamma_{\text{out}}$ are far away from the bifurcation and are not affected by $\mupar$.
Here $\u(\mupar)$ represents the unknown velocity, while $\p(\mupar)$ the unknown pressure. Moreover, the inlet velocity $\boldsymbol{u}_{\text{in}}$ is a parabolic profile and the viscosity $\nu$ is chosen such that the resulting Reynolds number is equal to 400, corresponding to the average Reynolds number over a cardiac cycle \cite{zarins1983carotid}.

As quantity of interest we would have liked to consider the pressure drop across the stenoses. However, due to the lack of physiological boundary conditions that prescribe the pressure at the inlet, as well as non-physiological homogeneous Neumann boundary conditions that prescribe zero pressure at the outlet, this quantity of interest would be severely affected by a variation of the degree of stenosis. To mitigate this, we resort to a relative pressure drop dividing by the inlet pressure (i.e., the pressure drop between the inlet and the outlet).
To be more precise, denote with $P_i(\mupar)$ the average of the pressure on the section $S_i$, $i = 0, \hdots, 5$, as shown in Figure \ref{fig:output_sections}. Sections $S_3$ and $S_4$ ($S_5$, respectively) are used to quantify the pressure drop for the stenosis on the left (right, resp.) branch, which is then divided by the pressure drop between the inlet $S_0$ and the left (right, resp.) outlet $S_1$ ($S_2$, resp.).
Therefore, the quantity of interest that we consider is the sum of the relative pressure drop of the two branches:
\begin{equation}
f(\mupar) = \frac{P_3(\mupar) - P_4(\mupar)}{P_0(\mupar) - P_1(\mupar)} + \frac{P_3(\mupar) - P_5(\mupar)}{P_0(\mupar) - P_2(\mupar)}.
\label{eq:qoi}
\end{equation}

\section{Shape morphing based on radial basis functions interpolation}
\label{sec:rbf}
Radial Basis Functions (RBF) are a powerful tool for shape
parametrization due to their good approximation
properties~\cite{buhmann2003radial, manzoni2012model}. 
In this section we summarize RBF-based shape morphing following the presentation in \cite{forti2014efficient}.
All the algorithms have been implemented
in the open source python package PyGeM~\cite{pygem}, which is used to perform the
shape morphing in the numerical results showed in Section~\ref{sec:results}.


A radial basis function is any smooth real-valued function $\widetilde{\varphi}: \mathbb{R}^n \to \mathbb{R}$
such that it exists $\varphi: \mathbb{R}^+ \to \mathbb{R}$ and
$\widetilde{\varphi} (\x) = \varphi (\| \x
\|)$, where $\| \cdot \|$ is the Euclidean norm in $\mathbb{R}^n$.

\begin{figure}[htb]
\centering
\includegraphics[width=0.6\textwidth]{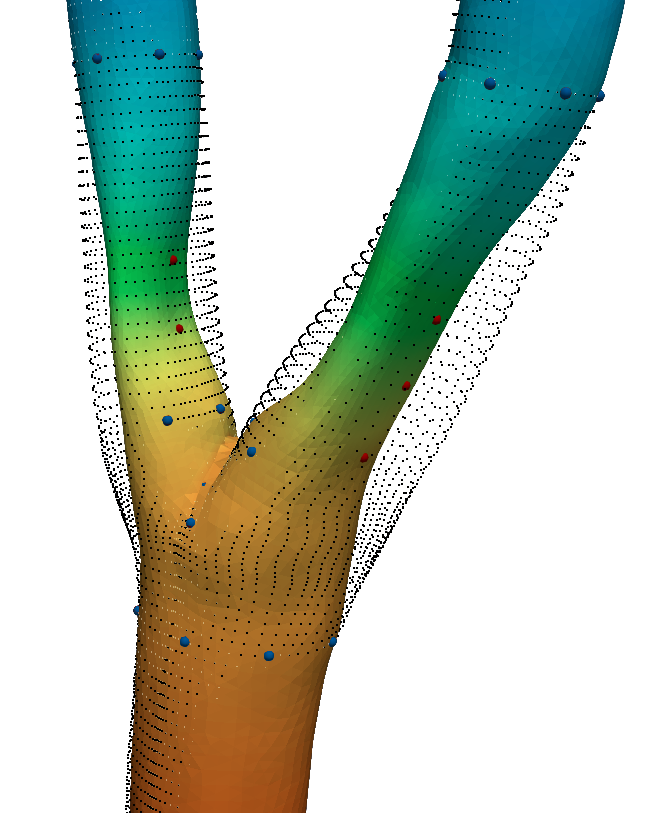}
\caption{Portion of a carotid. In black the original mesh represented
  by points, while in plain colors the deformed carotid.
  RBF interpolation control points are denoted by spherical markers (blue: fixed, red: not fixed).}
\label{fig:rbf_dormation}
\end{figure}

Let us recall that parameters are denoted by $\mupar \in \mathbb{D}$.
The RBF shape parametrization technique is based on the map
$\mathcal{M}(\x; \mupar) : \mathbb{R}^n \rightarrow \mathbb{R}^n$
, defined as follows 
\begin{equation}
\label{eq:rbf_map}
    \mathcal{M}(\x; \mupar) = q(\x; \mupar) + 
    \sum_{i=1}^{\mathcal{N}_C} \gamma_i(\mupar)\;
    \varphi(\| \x - \x_{C_i} \|).
\end{equation}
where $q(\x; \mupar)$ is a polynomial term to be determined,
$\{\gamma_i(\mupar)\}_{i=1}^{\mathcal{N}_C}$ are weights to be determined,
and $\{\x_{C_i}\}_{i=1}^{\mathcal{N}_C}$ are control points selected by the user
(denoted by spherical markers in Figure \ref{fig:rbf_dormation}), and
$\x \in \Omega$.
%
%
Among all the
possible RBF for modeling shapes we select the so-called thin plate
splines~\cite{duchon1977splines} defined as
\[
\varphi(r; R) =
            \left( \frac{r}{R} \right)^2
            \ln \left( \frac{r}{R} \right)
\]
where $r$ is the radial coordinate and $R > 0$ is a given radius.

In order to find the unknowns $q(\x; \mupar)$ and $\{\gamma_i(\mupar)\}_{i=1}^{\mathcal{N}_C}$, let
us assume that $q(\x; \mupar)$ is a polynomial function of degree 1, that is
\[
q(\x; \mupar) = \mathbf{c}(\mupar) + \mathbf{Q}(\mupar)\x,
\]
for some unknown $\mathbf{c}(\mupar) \in \mathbb{R}^n$ and $\mathbf{Q}(\mupar) \in \mathbb{R}^{n \times n}$. 
Therefore equation~\eqref{eq:rbf_map} can be rewritten in
matrix form as follows:
\begin{equation}
\label{eq:rbf_matrix_form}
    \mathcal{M}(\x; \mupar) = \mathbf{c}(\mupar) +
    \mathbf{Q}(\mupar)\x +
    \mathbf{G}^T(\mupar)\mathbf{d}(\x),
\end{equation}
being $\mathbf{d}(\x) = [\varphi(\| \x -
\x_{C_1} \|), \dots, \varphi(\| \x - \x_{\mathcal{N}_C}
\|)] \in \mathbb{R}^{\mathcal{N}_C}$ the vector constructed evaluating the
radial basis function on the Euclidean distance between the control
points position $\x_{C_i}$ and $\x$,
and the unknown $\mathbf{G}(\mupar) = [\gamma_1(\mupar),
\dots, \gamma_{\mathcal{N}_C}(\mupar)] \in \mathbb{R}^{\mathcal{N}_C \times
  n}$. 
To compute the unknowns $\mathbf{c}(\mupar)$, $\mathbf{Q}(\mupar)$ and $\mathbf{G}(\mupar)$
we enforce interpolation conditions on the set of control points, 
that is, given their
initial position as 
\[
\x_{C} = [\x_{C_1}, \dots,
\x_{C_{\mathcal{N}_C}}] \in \mathbb{R}^{\mathcal{N}_C \times n}
\]
and their $\mupar$-dependent deformed positions as
\begin{equation}
\mathbf{y}_{C}(\mupar) = [\mathbf{y}_{C_1}(\mupar), \dots,
\mathbf{y}_{C_{\mathcal{N}_C}}(\mupar)] \in \mathbb{R}^{\mathcal{N}_C \times n},
\label{eq:yC}
\end{equation}
we enforce that
\begin{equation}
\mathcal{M}(\x_{C_i}; \mupar) = \mathbf{y}_{C_i}(\mupar) \qquad \forall i \in
\{1, \dots, \mathcal{N}_C\}.
\label{eq:rbfls1}
\end{equation}
The system is then completed by additional
constraints that represent the conservation of the total force and
momentum~\cite{buhmann2003radial, dryden1998statistical, morris2008cfd}, due to
the presence of the polynomial term, as follows
\begin{align}
\sum_{i=1}^{\mathcal{N}_C} \gamma_i(\mupar) &= 0, \label{eq:rbfls2}\\
\sum_{i=1}^{\mathcal{N}_C} \gamma_i(\mupar) \x_{C_{1, i}} &= \dots =
\sum_{i=1}^{\mathcal{N}_C} \gamma_i(\mupar) \x_{C_{n, i}} = 0,\label{eq:rbfls3}
\end{align}
being $\x_{C_i} = [\x_{C_{1, i}}, \dots,
\x_{C_{n, i}}]$ a vector collecting the $i$-th coordinates
of all control points. These additional constraints, together with the presence of
the polynomial term of degree one, ensures that the resulting linear system~\eqref{eq:rbfls1}-\eqref{eq:rbfls3}
has always a unique solution ($\mathbf{c}(\mupar)$, $\mathbf{Q}(\mupar)$, $\mathbf{G}(\mupar)$).
Once the system is solved, we can deform all the points of the mesh through $\mathcal{M}(\mathbf{\cdot}; \mupar)$ to
obtain the deformed configuration.

\begin{figure}
\centering
\includegraphics[width=1.\textwidth]{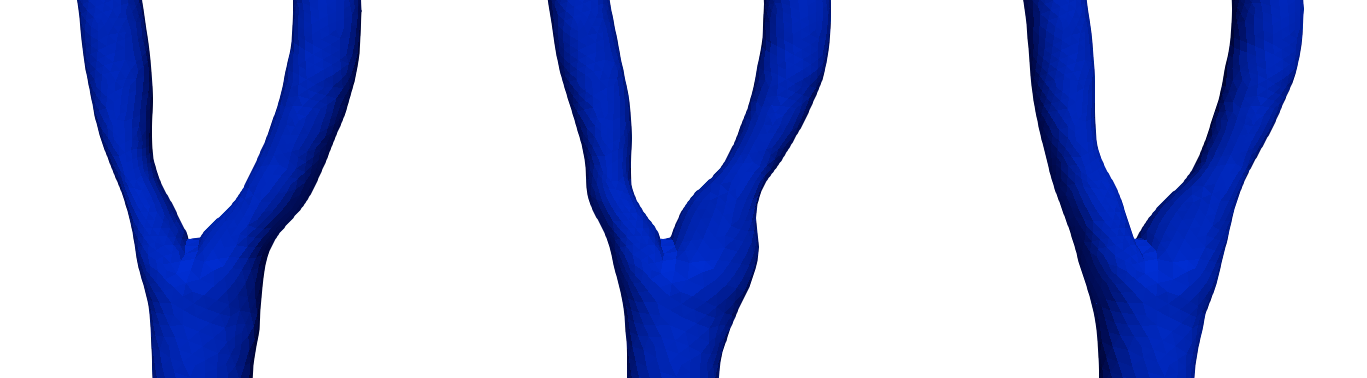}
\caption{Three different deformations produced varying the $m=10$ parameters.}
\label{fig:3dormations}
\end{figure}

In order to exemplify~\eqref{eq:yC} in our case, let us consider again Figure~\ref{fig:rbf_dormation}, where control
points are denoted by spherical markers. In order to enforce 
the deformation to be localized, we will move all control points close to the stenosis (colored in red) in the normal direction,
thus varying the occlusion,
while we will keep fixed a few control points far away from the stenosis (colored in blue). Therefore, for all fixed control points the
right hand side of \eqref{eq:rbfls1} is defined employing
\[
\mathbf{y}_{C_i}(\mupar) = \x_{C_i} \qquad \forall i \in \{1, \dots, \mathcal{N}_C^{\text{fixed}}\},
\]
while for all non-fixed control points
\[
\mathbf{y}_{C_i}(\mupar) = \x_{C_i} - \mu_i \mathbf{n}_i \qquad \forall i \in \{1, \dots, m\},
\]
where $\mathcal{N}_C = \mathcal{N}_C^{\text{fixed}} + m$, $\mu_i$ denotes the $i$-th element in $\mupar$ and $\mathbf{n}_i$
the outer unit normal to the wall evaluated at $\x_{C_i}$. 
In our case we have chosen $m=10$, $\mathcal{N}_C^{\text{fixed}}=55$ (not all control points are shown in Figure~\ref{fig:rbf_dormation}),
and the parameter range is $\mathbb{D} = [0, 0.3]^m$, resulting in a
wide range of possible different stenosis scenarios (see
Figure~\ref{fig:3dormations} for an idea of the possible configurations).

\section{High fidelity solver based on the finite element method}
\label{sec:solver} 
In this section we summarize the high fidelity solver that will be
used in the training of both the active subspace (see Section
\ref{sec:as}) and the POD-Galerkin reduction (see Section
\ref{sec:pod}). Let $(V, Q)$ be an inf-sup stable finite element (FE)
pair, being $V_{\delta}(\mupar) \subset
H^1(\Omega_{\delta}(\mupar);\;\mathbb{R}^n)$ and $Q_{\delta}(\mupar)
\subset L^2(\Omega_{\delta}(\mupar))$, being $\delta$ the maximum
diameter size of a tetrahedralization\footnote{In order to simplify the exposition we will report the FE formulation on the deformed domain $\Omega(\mupar)$. However, it should be noted that only the mesh $\Omega_\delta$ of the reference domain $\Omega$ is generated, and deformed meshes $\Omega_\delta(\mupar)$ are obtained through the mapping $\mathcal{M}(\boldsymbol{\cdot}; \mupar)$.} $\Omega_{\delta}$ of $\Omega$; Taylor-Hood $\mathbb{P}^2-\mathbb{P}^1$ FE has been employed in the numerical experiments shown in Section \ref{sec:results}.

\begin{figure}[h!]
\centering
\subfloat[Velocity field]{\includegraphics[height=0.6\textwidth]{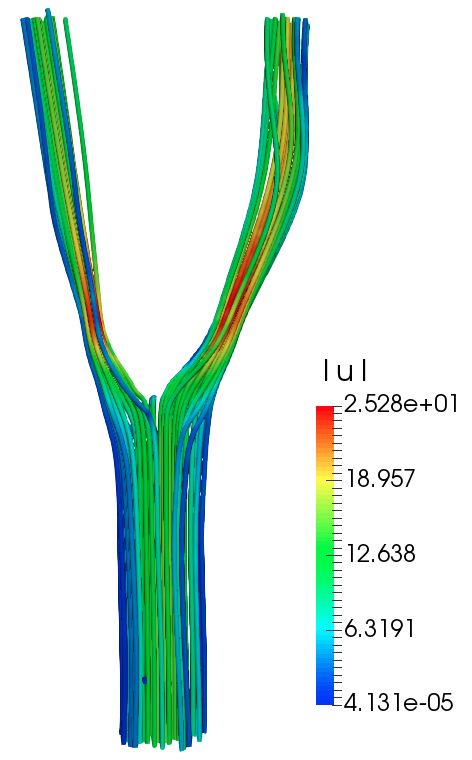}\label{subfig:corner_vel}}\hspace{.1cm}
\subfloat[Vorticity]{\includegraphics[height=0.6\textwidth]{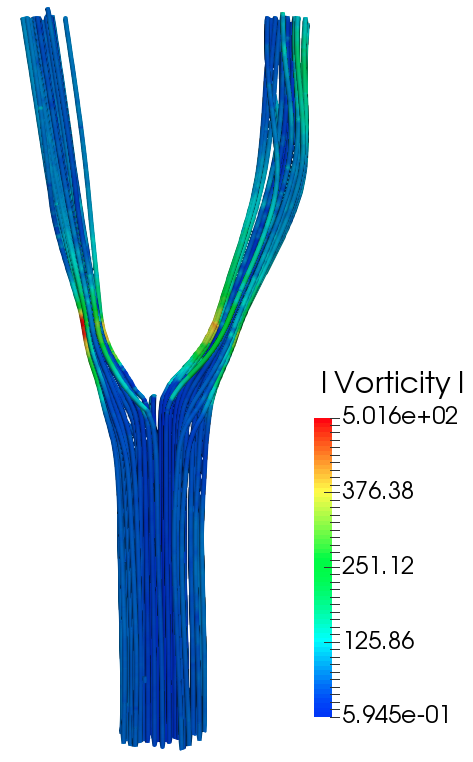}\label{subfig:corner_vor}}
\caption{Velocity and vorticity fields for a representative choice of $\mupar \in \mathbb{D}$.}
\label{fig:corner_case}
\end{figure}

The high fidelity approximation of \eqref{eq:ns} reads: for any $\mupar~\in~\mathbb{D}$, find $(\u_{\delta}(\mupar), \p_{\delta}(\mupar)) \in V_{\delta}(\mupar) \times Q_{\delta}(\mupar)$ such that $\u_\delta(\mupar) = \boldsymbol{u}_{\text{in}} \text{ on } \Gamma_{\text{in}}$, $\u_\delta(\mupar) = \boldsymbol{0}, \text{ on } \Gamma_{\text{wall}}(\mupar)$ and:
\begin{equation}\label{eq:nsdelta}
\begin{cases}
a(\u_\delta, \boldsymbol{v}_\delta; \mupar) + b( \boldsymbol{v}_\delta, \p_\delta; \mupar ) + c(\u_\delta, \u_\delta,  \boldsymbol{v}_\delta; \mupar) = \boldsymbol{0}, & \forall  \boldsymbol{v}_\delta \in V_\delta(\mupar),\\
b(\u_\delta, q_\delta; \mupar) = 0, & \forall q \in Q_\delta(\mupar),
\end{cases}
\end{equation}
where
\begin{equation*}
a(\u_\delta, \boldsymbol{v}_\delta ;\mupar)=  \int_{\Omega(\mupar)} \nu \nabla\u_\delta : \nabla \boldsymbol{v}_\delta \, d{\x}, \qquad 
b(\boldsymbol{v}_\delta, q_\delta; \mupar)= - \int_{\Omega(\mupar)} q_\delta \; \text{div}\; \boldsymbol{v}_\delta\, d{\x}
 \end{equation*}
are the bilinear forms associated to the diffusion and divergence operators, respectively, whereas
\begin{equation*}
c(\u_\delta, \boldsymbol{v}_\delta, \boldsymbol{z}_\delta; \mupar) =  \int_{\Omega(\mupar)} 
(\nabla\boldsymbol{v}_\delta \; \u_\delta) \cdot \boldsymbol{z}_\delta \, d{\x}
\end{equation*}
is the trilinear form related to the nonlinear advection term. Let $\{\boldsymbol{\varphi}^\delta_i\}_{i=1}^{N^\delta_{\bf u}}$ and $\{\zeta^\delta_k\}_{k=1}^{N^\delta_{p}}$  the bases of $V_\delta(\mupar)$ and $Q_\delta(\mupar)$, respectively. Then, \eqref{eq:nsdelta} can be equivalently rewritten as the following nonlinear system: for any $\mupar~\in~\mathbb{D}$, find $(\mathbf{u}(\mupar), \mathbf{p}(\mupar)) \in \mathbb{R}^{N^\delta_{\bf u}} \times \mathbb{R}^{N^\delta_{p}}$ such that\footnote{The non-homogeneous right-hand side accounts for boundary conditions via a lifting.}:
\begin{equation}\label{eq:nsboldsymbol}
\left[
\begin{array}{cc}
 \mathbf{A}(\mupar)  +    \mathbf{C}({{\bf u}}(\mupar); \mupar) &  \mathbf{B}^T(\mupar) \\
 \mathbf{B}(\mupar) & \mathbf{0}
\end{array}
\right]
\left[
\begin{array}{c}
{{\bf u}}(\mupar) \\ {{\bf p}}(\mupar)\end{array}
\right]
=
\left[
\begin{array}{c}
{{\bf f}}(\mupar) \\ {{\bf g}}(\mupar)
\end{array}
\right],
\end{equation}
where, for $1 \leq i,j \leq N^\delta_{\bf u}$ and $1 \leq k \leq N^\delta_{p}$:
\begin{equation}\label{eq:nsboldsymboldef}
\begin{array}{c}
  (\mathbf{A}(\mupar))_{ij} = a(\boldsymbol{\varphi}^\delta_j, \boldsymbol{\varphi}^\delta_i; \mupar), \qquad (\mathbf{B}(\mupar))_{ki} = b(\boldsymbol{\varphi}^\delta_i, \zeta^\delta_k; \mupar), \vspace{0.25cm}\\
(\mathbf{C}({\bf u} (\mupar); \mupar))_{ij} = c\left(\boldsymbol{u}_\delta(\mupar), \boldsymbol{\varphi}^\delta_j, \boldsymbol{\varphi}^\delta_i; \mupar\right),
\end{array}
\end{equation}
and
\begin{equation*}
\boldsymbol{u}_\delta(\mupar) = \sum_{i=1}^{N^\delta_{\bf u}} u_i(\mupar) \boldsymbol{\varphi}^\delta_i, \quad
p_\delta(\mupar) = \sum_{k=1}^{N^\delta_{\bf p}} p_k(\mupar) \zeta^\delta_k,
\end{equation*}
For any $\mupar \in \mathbb{D}$ the nonlinear system \eqref{eq:nsboldsymbol} is solved via a Newton method, and the quantity of interest $f(\mupar)$ is computed according to \eqref{eq:qoi} in a postprocessing stage. The solution of the Navier-Stokes problem for a representative value of $\mupar$ is shown in Figure \ref{fig:corner_case}.

\section{Parameter space reduction by active subspaces property}
\label{sec:as}
The active subspaces (AS) property has been emphasized recently, see
for example P. Constantine in~\cite{constantine2015active}. It concerns the properties
of a parametric scalar function and it is exploited for
dimension reduction in parameter studies. 
The main idea behind the active subspaces is the following: we rescale
the inputs (the parameters $\mupar$, in our case) and then rotate the inputs domain with respect to the
origin in such a way lower dimension behavior of the output function (the normalized pressure drop $f(\mupar)$, in our case)
is revealed. The active subspaces identify a set of important
directions in the space of all inputs, instead of identifying a
subset of the inputs as important. The latter approach would have indeed resulted in
similar limitations as in our previous works~\cite{BallarinManzoniRozzaSalsa2013,forti2014efficient,DAmario}.
If the output of the simulation
does not change on average along a particular direction of the
parameters, then we can safely ignore that direction in the parameter
study. When an active
subspace is identified for the problem of interest, it is
possible to perform different parameter studies such as integration, optimization, response
surfaces~\cite{box1987empirical}, and statistical inversion~\cite{kaipio2006statistical}.

Now we review the process of finding active subspaces. Let us assume\footnote{In this section we will omit the dependence on $\mupar$. It should be understood that $f = f(\mupar)$, $\rho = \rho(\mupar)$, etc.}
$f: \mathbb{R}^m \rightarrow \mathbb{R}$ is a scalar function and $\rho: \mathbb{R}^m \rightarrow \mathbb{R}^+$ a probability density function.
For our sake a uniform probability density will suffice, as all possible geometrical configurations can be drawn with equal probability.
In particular, we assume $f$ continuous and differentiable in the support of $\rho$, with
continuous and square-integrable (with respect to the measure induced by $\rho$) derivatives. The
active subspaces of the pair $(f, \rho)$ are the eigenspaces of the
covariance matrix associated to the gradients $\nabla_{\mupar} f$. 
To this end we define the
so-called uncentered covariance matrix of the gradients of $f$ (among others see~\cite{devore2015probability} for a
more deep understanding of these operators), denoted by $\boldsymbol{\Sigma}$, whose elements are the average products of partial
derivatives of the simulations' input/output map, i.e.:
\begin{equation}
\label{eq:covariance}
\boldsymbol{\Sigma} = \mathbb{E}\, [\nabla_{\mupar} f \, \nabla_{\mupar} f
^T] = \int_{\mathbb{D}} (\nabla_{\mupar} f) ( \nabla_{\mupar} f )^T
\rho \, d \mupar ,
\end{equation}
where $\mathbb{E}[\cdot]$ is the expected value.
To approximate the eigenpairs of this matrix it is common to use a Monte Carlo method as follows \cite{metropolis1949monte,constantine2015computing}:
\begin{equation}
\label{eq:covariance_approx}
\boldsymbol{\Sigma} \approx \frac{1}{N_{\text{train}}^{\text{AS}}} \sum_{i=1}^{N_{\text{train}}^{\text{AS}}} \nabla_{\mupar} f_i \,
\nabla_{\mupar} f^T_i ,
\end{equation}
where we draw $N_{\text{train}}^{\text{AS}}$ independent samples $\mupar^{(i)}$ from the
measure $\rho$ and where $\nabla_{\mupar} f_i = \nabla_{\mupar}
f(\mupar^{(i)})$.
The matrix $\boldsymbol{\Sigma}$
is symmetric and positive semidefinite, so it admits a real eigenvalue
decomposition
\begin{equation}
\label{eq:decomposition}
\boldsymbol{\Sigma} = \mathbf{W} \mathbf{\Lambda} \mathbf{W}^T ,
\end{equation}
where $\mathbf{W}$ is a $m \times m$ orthogonal matrix of eigenvectors,
and $\mathbf{\Lambda}$ is the diagonal matrix of non-negative
eigenvalues arranged in descending order.

We now form a lower dimensional parameter subspace by selecting the first $M$ eigenvectors, for some
$M < m$. 
On average, perturbations
in the first set of coordinates change $f$ more than perturbations in
the second set of coordinates. While low eigenvalues suggest
that the corresponding vector is in the nullspace of the covariance
matrix. Discarding these vectors we can construct an approximation of
$f$. For the sake of notation, let us partition $\mathbf{\Lambda}$ and $\mathbf{W}$ as follows:
\[
\mathbf{\Lambda} =   \begin{bmatrix} \mathbf{\Lambda}_1 & \\
                                     &
                                     \mathbf{\Lambda}_2\end{bmatrix},
\qquad
\mathbf{W} = \left [ \mathbf{W}_1 \quad \mathbf{W}_2 \right ],
\]
where $\mathbf{\Lambda}_1 = \text{diag}(\lambda_1, \dots, \lambda_M)$, and
$\mathbf{W}_1$ contains the first $M$ eigenvectors. The active subspace
is the the range of $\mathbf{W}_1$. 
The inactive subspace is the range of the remaining eigenvectors in
$\mathbf{W}_2$.
The linear combinations of the input
parameters with weights from the important eigenvectors are the active variables. 
We approximate the behaviour of the objective
function by projecting the full parameter space onto the active
subspace. 

Given the input parameters in the full space we can map forward to the
active subspace. Respectively we have the following formulas for the
active variable $\mupar_M$ and the inactive variable $\etapar$:
\begin{equation}
\label{eq:active_var}
\mupar_M = \mathbf{W}_1^T\mupar \in \mathbb{R}^M, \qquad
\etapar = \mathbf{W}_2^T \mupar \in \mathbb{R}^{m - M} .
\end{equation}
That means that any point in the parameter space $\mupar \in
\mathbb{R}^m$ can be expressed in terms of $\mupar_M$ and
$\etapar$ as follows:
\[
\mupar = \mathbf{W}\mathbf{W}^T\mupar =
\mathbf{W}_1\mathbf{W}_1^T\mupar +
\mathbf{W}_2\mathbf{W}_2^T\mupar = \mathbf{W}_1 \mupar_M +
\mathbf{W}_2 \etapar.
\]
So it is possible to rewrite $f$ as
\[ 
f (\mupar) =  f (\mathbf{W}_1 \mupar_M + \mathbf{W}_2 \etapar) ,
\]
and construct a surrogate quantity of interest $g$ using only the active variables
\[
f (\mupar) \approx g (\mathbf{W}_1^T \mupar) = g(\mupar_M).
\]
In our pipeline, the surrogate quantity of interest $g$ will be obtained querying a POD-Galerkin reduced
order model, as described in the next section.

Active subspaces can be seen in the more general context of ridge
approximation (see for example~\cite{pinkus2015ridge,keiper2015analysis}). It can be proved that, under certain
conditions, the active subspace is a good starting point in optimal
ridge approximation and it is nearly stationary as shown
in~\cite{constantine2016near,hokanson2017data}.

\section{Model order reduction based on a POD-Galerkin method}
\label{sec:pod}
In this section we now briefly summarize the POD-Galerkin method that we employ for the model order reduction of the high fidelity approximation~\eqref{eq:nsboldsymbol} of the Navier-Stokes equations \eqref{eq:ns}, based on the usual offline-online paradigm~\cite{MSA9,MSA17,BallarinManzoniQuarteroniRozza2014,hesthaven2015certified}. The main novelty in this section is that the training during the
offline stage (as well as the testing during the online one) will be carried out \emph{only} over the active (parameter) subspace, and not
over the full parameter space~$\mathbb{D}$.

Let us denote by $\Xi_{\text{train}}^{\text{POD}} = \{\mupar_M^{(i)}\}_{i=1}^{N_{\text{train}}^{\text{POD}}} \subset \mathbb{D}$ a training set of $N_{\text{train}}^{\text{POD}}$ points chosen randomly over the active subspace, i.e. the range of $\mathbf{W}_1$. 
During the offline stage, we assemble the following \emph{snapshots} matrices:
\begin{align*}
&S_{\bf u} = [{\bf u}(\mupar_M^{(1)}) \ | \  \ldots \ | \  {\bf u}(\mupar_M^{(N_{\text{train}})}) ] \in \mathbb{R}^{N_{\bf u}^\delta \times N_{\text{train}}^{\text{POD}}},\\
&S_{p} = [ {\bf p}(\mupar_M^{(1)}) \ | \ \ldots \ | \ {\bf p}(\mupar_M^{(N_{\text{train}})}) ] \in \mathbb{R}^{N_{p}^\delta \times N_{\text{train}}^{\text{POD}}},\\
&S_{\bf s} = [{\bf s}(\mupar_M^{(1)}) \ | \  \ldots \ | \  {\bf s}(\mupar_M^{(N_{\text{train}})}) ] \in \mathbb{R}^{N_{\bf u}^\delta \times N_{\text{train}}^{\text{POD}}},
\end{align*}
where $({\bf u}(\boldsymbol{\cdot}), {\bf p}(\boldsymbol{\cdot}))$ is the FE solution of \eqref{eq:nsboldsymbol}, and the supremizer solution ${\bf s}(\boldsymbol{\cdot})$ is obtained by a FE approximation of the following elliptic equation: for each $\mupar \in \Xi_{\text{train}}^{\text{POD}}$, assuming (an approximation of) $\p_\delta(\mupar)$ to be known, find $\boldsymbol{s}_\delta(\mupar) \in V_\delta(\mupar)$ such that
\begin{equation*}
(\boldsymbol{s}_\delta(\mupar), \boldsymbol{v}_\delta)_{V_\delta(\mupar)} =  b( \boldsymbol{v}_\delta, \p_\delta(\mupar); \mupar ) \quad \forall  \boldsymbol{v}_\delta \in V_\delta(\mupar),
\end{equation*}
where $(\boldsymbol{\cdot}, \boldsymbol{\cdot})_{V_\delta(\mupar)}$ represents the inner products in $H^1(\Omega(\mupar);\;\mathbb{R}^n)$.
Indeed, $\boldsymbol{s}_\delta(\mupar)$ is such that
\begin{equation*}
\boldsymbol{s}_\delta(\mupar) = \arg \sup_{\boldsymbol{v}_\delta \neq \boldsymbol{0}} \frac{b( \boldsymbol{v}_\delta, \p_\delta(\mupar); \mupar ) }{\|\boldsymbol{v}_\delta\|_{V_\delta(\mupar)}},
\end{equation*}
and the inf-sup constant
\begin{equation*}
{\beta}_\delta(\mupar)  = \inf_{q_\delta \neq 0} \sup_{\boldsymbol{v}_\delta \neq \boldsymbol{0}} \frac{b( \boldsymbol{v}_\delta, \p_\delta(\mupar); \mupar ) }{\|\boldsymbol{v}_\delta\|_{V_\delta(\mupar)} \|q_\delta\|_{Q_\delta(\mupar)}}
\end{equation*}
is related to the the supremizer solution as follows:
\begin{equation*}
\left({\beta}_\delta(\mupar)\right)^2 = \inf_{q_\delta \neq 0} \frac{\|\boldsymbol{s}_\delta(\mupar)\|_{V_\delta(\mupar)}}{\|q_\delta\|_{Q_\delta(\mupar)}}.
\end{equation*}
Such supremizers are employed at the reduced order level to enhance the inf-sup stability of the reduced system, which is essential to obtain an accurate approximation of the pressure. We refer to \cite{RozzaVeroy2007,BallarinManzoniQuarteroniRozza2014} for further details on supremizers, as well as to \cite{CaiazzoIliescuJohnSchyschlowa2013,Stabile2018,Shafqat} for possible alternative approaches.

A POD basis for the velocity, pressure and supremizer spaces can be obtained by considering the singular value decomposition of the following matrices \cite{ravindran00,Berkooz_POD_review}
\[
X_{\bf u}^{1/2} S_{\bf u}, \qquad X_{p}^{1/2} S_{p}, \qquad X_{\bf u}^{1/2} S_{\bf s},
\]
being $X_{\bf u}$ and $X_{p}$ FE matrices corresponding to the discretization of the inner products in $H^1(\Omega;\;\mathbb{R}^n)$ and $L^2(\Omega)$, respectively.
The first $N_{\bf u}, N_p, N_{\bf s}$ (respectively) left singular vectors are then considered as basis functions $\{\boldsymbol{\varphi}_n\}_{n=1}^{N_{\bf u}}, \{\zeta_n\}_{n=1}^{N_{p}}, \{\boldsymbol{\phi}_n\}_{n=1}^{N_{\bf s}}$ of the reduced basis spaces. Therefore, the reduced spaces for velocity $V_N$ and pressure $Q_N$, of cardinality $N_{{\bf u},{\bf s}} = N_{\bf u} + N_{\bf s}$ and $N_p$, respectively,
are then obtained as 
\[
V_N = \text{span}(\{\boldsymbol{\varphi}_n\}_{n=1}^{N_{\bf u}}, \{\boldsymbol{\phi}_n\}_{n=1}^{N_{\bf s}}), \quad
Q_N = \text{span}(\{\zeta_n\}_{n=1}^{N_{p}}).
\] 
Finally, let us introduce the corresponding basis functions matrices
\begin{align*}
& \mathbf{Z}_{{\bf u},{\bf s}} = [\boldsymbol{\varphi}_1 | \hdots | \boldsymbol{\varphi}_{N_{\bf u}} | \boldsymbol{\phi}_1 | \hdots | \boldsymbol{\phi}_{N_{\bf s}}] \in \mathbb{R}^{N_{\bf u}^\delta \times N_{{\bf u},{\bf s}}},\\
& \mathbf{Z}_{p} = [\zeta_1 | \hdots | \zeta_{N_{p}}] \in \mathbb{R}^{N_{p}^\delta \times N_{p}}.
\end{align*}

Let us now denote by $\Xi_{\text{test}}^{\text{POD}} = \{\mupar_M^{(j)}\}_{j=1}^{N_{\text{test}}^{\text{POD}}} \subset \mathbb{D}$ a testing set of $N_{\text{test}}^{\text{POD}}$ points chosen randomly over the active subspace. During the online stage, for any $\mupar_M \in \Xi_{\text{test}}^{\text{POD}}$, we solve the following reduced nonlinear system: find $(\mathbf{u}_N(\mupar_M), \mathbf{p}_N(\mupar_M)) \in \mathbb{R}^{N_{\bf u, \bf s}} \times \mathbb{R}^{N_{p}}$ such that
\begin{equation*}
\left[
\begin{array}{cc}
 \mathbf{A}_N(\mupar_M)  +    \mathbf{C}_N({{\bf u}}_N(\mupar_M); \mupar_M) &  \mathbf{B}_N^T(\mupar_M) \\
 \mathbf{B}_N(\mupar_M) & 0
\end{array}
\right]
\left[
\begin{array}{c}
{{\bf u}}_N(\mupar_M) \\ {{\bf p}}_N(\mupar_M)\end{array}
\right]
=
\left[
\begin{array}{c}
{{\bf f}}_N(\mupar_M) \\ {{\bf g}}_N(\mupar_M)
\end{array} \vspace{-0.15cm}
\right],
\end{equation*}
where
\begin{align*}
&\mathbf{A}_N(\mupar_M) = \mathbf{Z}_{\bf u, \bf s}^T \; \mathbf{A}(\mupar_M) \; \mathbf{Z}_{\bf u, \bf s}, \quad 
\mathbf{B}_N(\mupar_M)  = \mathbf{Z}_{p}^T \; \mathbf{B}(\mupar_M) \; \mathbf{Z}_{\bf u, \bf s}, \\ 
&\mathbf{C}_N({\bf w} ; \mupar_M) = Z_{\bf u, \bf s}^T \; \mathbf{C}(\mathbf{Z}_{\bf u, \bf s} \; {\bf w}; \mupar_M) \; \mathbf{Z}_{\bf u, \bf s}.
\end{align*}
In order to obtain the maximum efficiency, the reduced system should not require evaluations of quantities defined on the high fidelity mesh; a standard approach based on the empirical interpolation method \cite{BarraultMadayNguyenPatera2004} and its discrete variants, as well as alternatives based on gappy POD \cite{CarlbergFarhatCortialAmsallem2013}, could be considered. We omit any additional detail on this topic for the sake of brevity, as this is already extensively discussed in the literature cited in this section.

\section{Numerical results}
\label{sec:results}
In this section we present the results of the complete pipeline applied
to a specific artery bifurcation. Moreover we demonstrate the improvements
obtained using the pipeline with respect to the POD approach on the
full parameter space.

The mesh is discretized using tetrahedral cells; a FE approximation by $\mathbb{P}^2-\mathbb{P}^1$ elements is used,
resulting in 265049 degrees of freedom. FEniCS is employed for the implementation of the 
high fidelity solver \cite{LoggMardalEtAl2012a}.

Let us recall that the parameter space is a $m=10$ dimensional space. In particular the
parameters are the displacements of 10 different RBF control points
along the orthogonal direction with respect to the surface. The moving
control points are located in the two branches just after the bifurcation in
order to simulate a stenosed carotid. In Figure
\ref{fig:rbf_dormation} it is possible to observe the original
undeformed carotid in black and the moving control points in red with a
possible deformation. The PyGeM open source package is used to perform
the deformation \cite{pygem}. 

\begin{figure}[ht]
\centering
\subfloat[Maximum, minimum, and mean value of the aspect ratio for all
the mesh deformations.]{\includegraphics[width=0.475\textwidth]{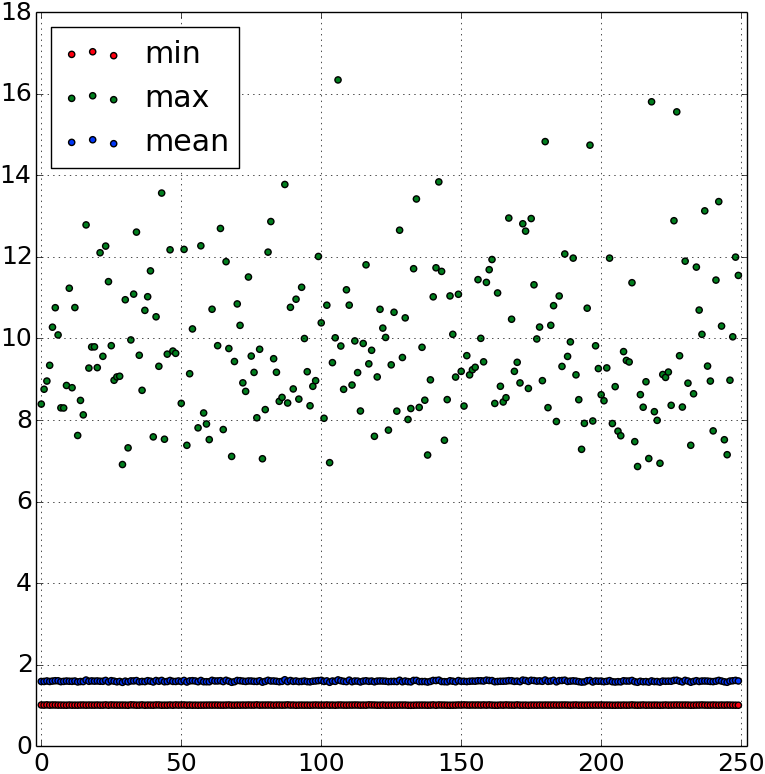}\label{subfig:aspect}}\hfill
\subfloat[Percentage of tetrahedra with an aspect ratio above the
maximum value of the reference geometry.]{\includegraphics[width=0.505\textwidth]{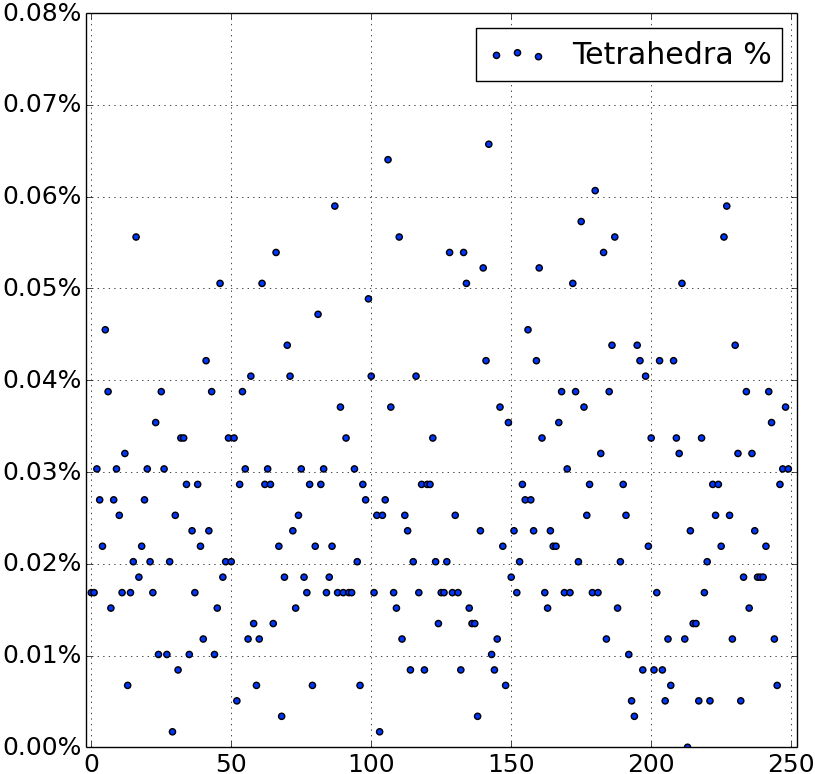}\label{subfig:aspect_above}}
\caption{Aspect ratio for each deformed mesh in the training set
  \protect\subref{subfig:aspect} and percentage of cells above the
  maximum aspect ratio of the reference mesh with respect to the total
  number of cells \protect\subref{subfig:aspect_above}.}
\label{fig:aspect_ratio}
\end{figure}

Since the deformations are made with respect to the reference
geometry, the quality of the resulting mesh after interpolation could
decrease. To address this problem we computed the aspect
ratio~\cite{frey2000mesh} of all the tetrahedra of each deformed
mesh. In Figure~\ref{subfig:aspect} we plot the minimum, the maximum
and the mean of such ratio. Even though the maximum values of such index are increased by the deformation process,
a sensible deterioration in mesh quality affects at most $0.07\%$ of the total number of
tetrahedra. Results are reported in Figure~\ref{subfig:aspect_above}, which summarizes the percentage of cells
for which the aspect ratio is above the maximum value it had in the reference configuration.
Thus, we conclude that the deformations impact on the mesh quality is
negligible.

Recalling Section~\ref{sec:as}, we need to construct the uncentered
covariance matrix $\boldsymbol{\Sigma}$ defined in \eqref{eq:covariance}. As shown in
\eqref{eq:covariance_approx}, we use a Monte Carlo method,
in order to construct the matrix $\boldsymbol{\Sigma}$, using the software \cite{paul-contantine-16}.
The number of training samples that we employ is $N_{\text{train}}^{\text{AS}} = 250$.
Even though it may be challenging to explore a 10 dimensional space, heuristics reported in~\cite{constantine2015active} suggest
this choice of $N_{\text{train}}^{\text{AS}}$ is enough for the purposes of the active subspaces identification.
In order to
approximate the gradients of the pressure drop $f$ with respect to the
parameters, that is $\nabla_{\mupar} f$, we use a local linear
model that approximates the gradients with the best linear
approximation using 17 nearest neighbors.
After constructing the matrix $\mathbf{\Sigma}$ we calculate its real eigenvalue decomposition.

\begin{figure}
\centering
\includegraphics[width=0.65\textwidth]{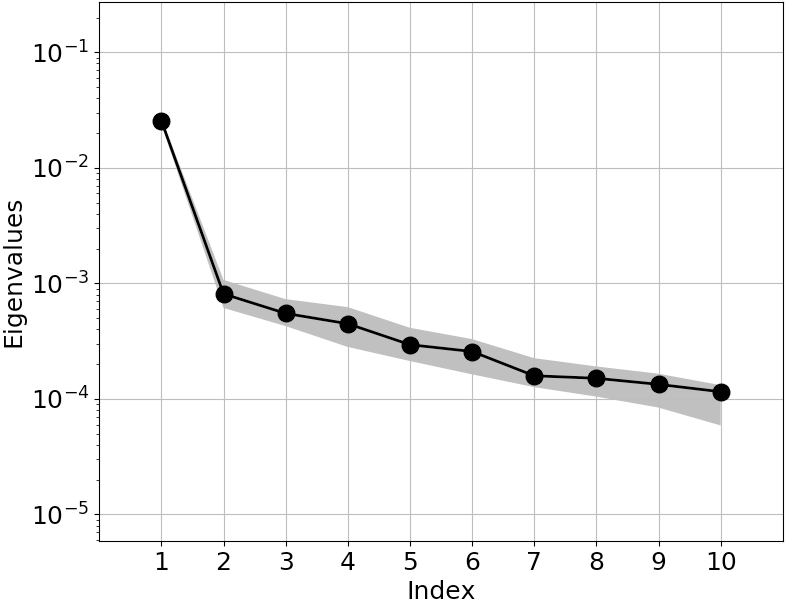}
\caption{Eigenvalue estimates in block circles with the bootstrap
  intervals (grey region). The order-of-magnitude gaps between the
  eigenvalues suggest confidence in the dominance of the active
  subspace.}
\label{fig:eigs_as}
\end{figure}

\begin{figure}[ht]
\centering
\subfloat[One active variable.]{\includegraphics[width=0.43\textwidth]{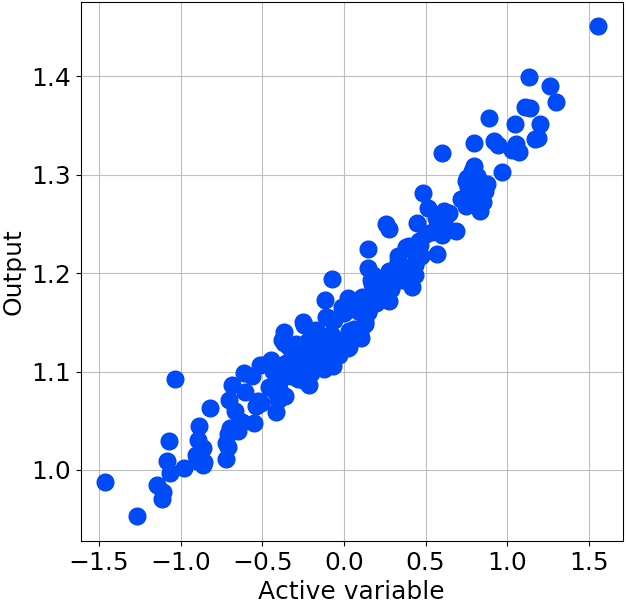}\label{subfig:ssp1}}\hfill
\subfloat[Two active variables.]{\includegraphics[width=0.523\textwidth]{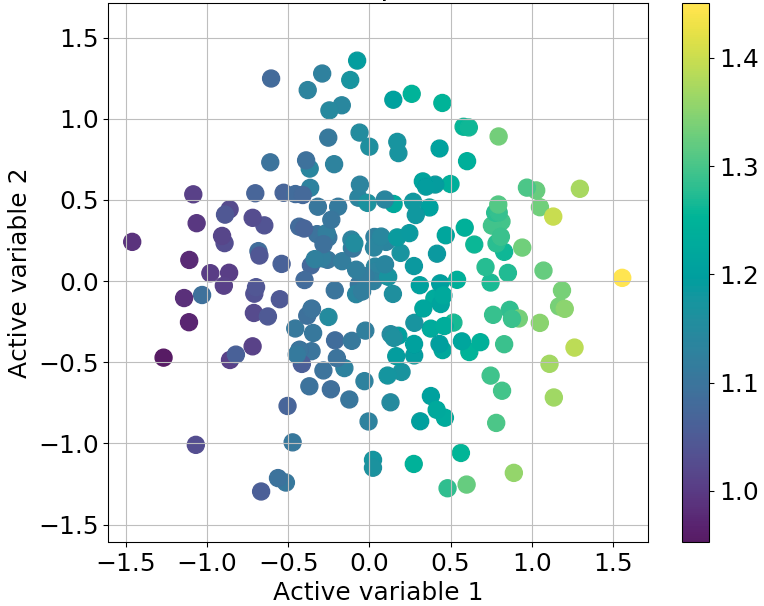}\label{subfig:ssp2}}
\caption{Sufficient summary plots for \protect\subref{subfig:ssp1} one and \protect\subref{subfig:ssp2} two active variables using the training dataset.}
\label{fig:ssp12}
\end{figure}

The eigenvalues of the covariance matrix in descending order are depicted in
Figure~\ref{fig:eigs_as}. The presence of gaps between the eigenvalues
supports the existence of an active subspace.
We can investigate the proper dimension of the active
subspace using scatter plots that contain all available regression
information that are called sufficient summary
plots~\cite{cook2009regression}. Recalling~\eqref{eq:active_var}, Figure~\ref{fig:ssp12} shows
$f (\mupar)$ against $\mupar_M = \mathbf{W}_1^T\mupar$, where
$\mathbf{W}_1^T$ contains the first one and the first two eigenvectors,
respectively. An active subspace of dimension one could suffice, but
the band-width of the scatter points is quite large, so we prefer to
retain more information about the output function by using a two
dimensional active subspace.

\begin{figure}
\centering
\includegraphics[width=0.6\textwidth]{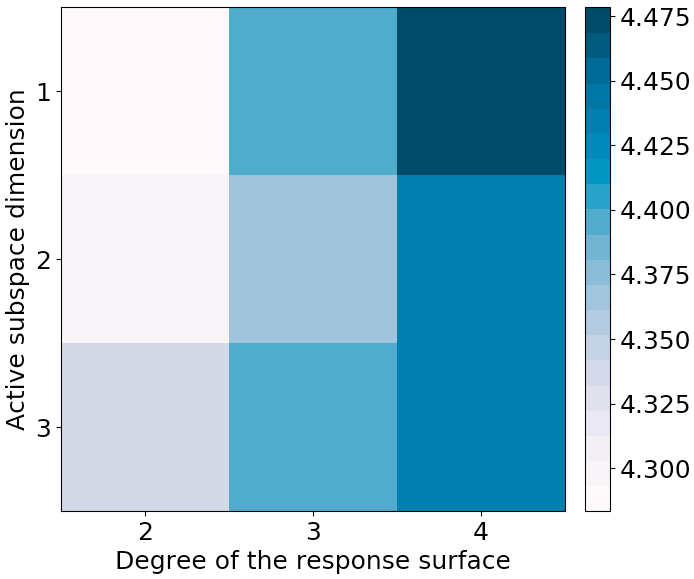}
\caption{Surrogate model error with respect to the active subspace dimension and the order of the response surface.}
\label{fig:heatmap}
\end{figure}

To support this decision we construct a response surface varying the
dimension of the active subspace and the order of the polynomial
surface and we compute the relative error with respect to a test
dataset. We can see from Figure~\ref{fig:heatmap} that on average the
bidimensional subspace is the best choice in terms of information
retention and dimension of the reduced parameter space. We underline
that the choice of the active subspace dimension depends on the
problem, the accuracy, and the goal you want to achieve. For the purpose of this
chapter the choice we made is a very good compromise and does not
affect the following results.

\begin{figure}[H]
\centering
\subfloat[Velocity.]{\includegraphics[width=0.48\textwidth]{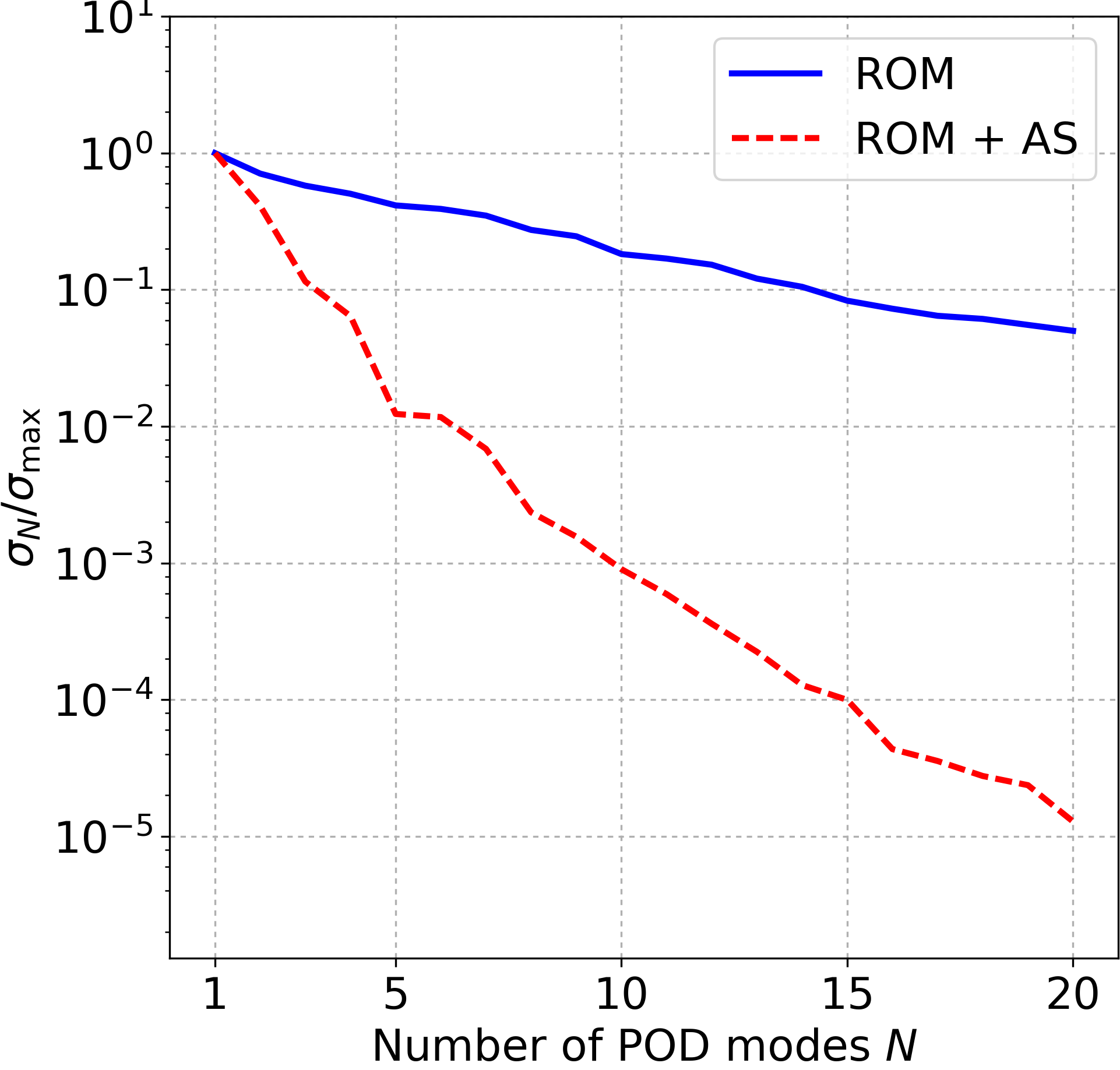}\label{subfig:vel_sval}}
\hfill
\subfloat[Supremizers.]{\includegraphics[width=0.48\textwidth]{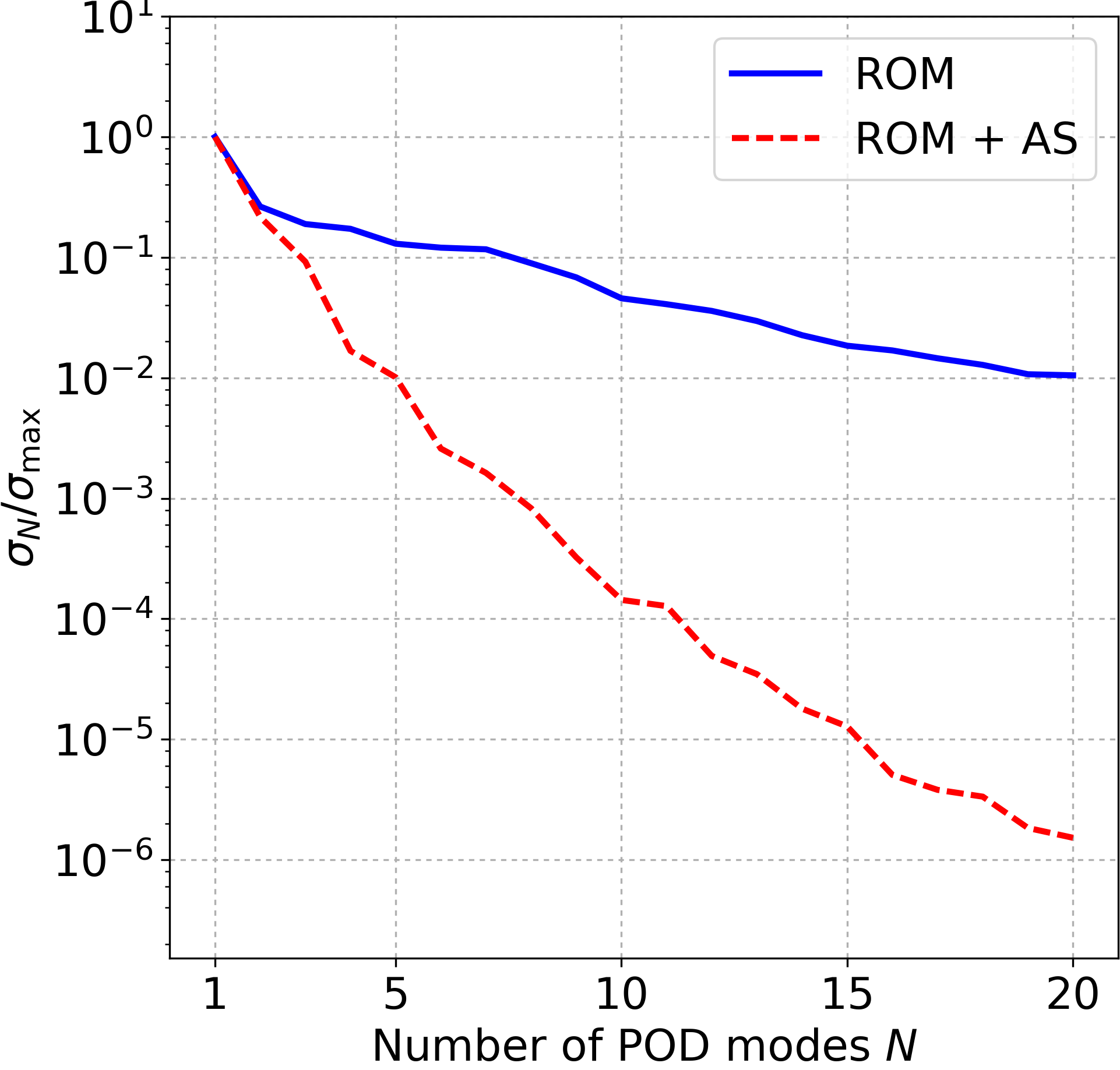}\label{subfig:sup_sval}}
\\
\subfloat[Pressure.]{\includegraphics[width=0.48\textwidth]{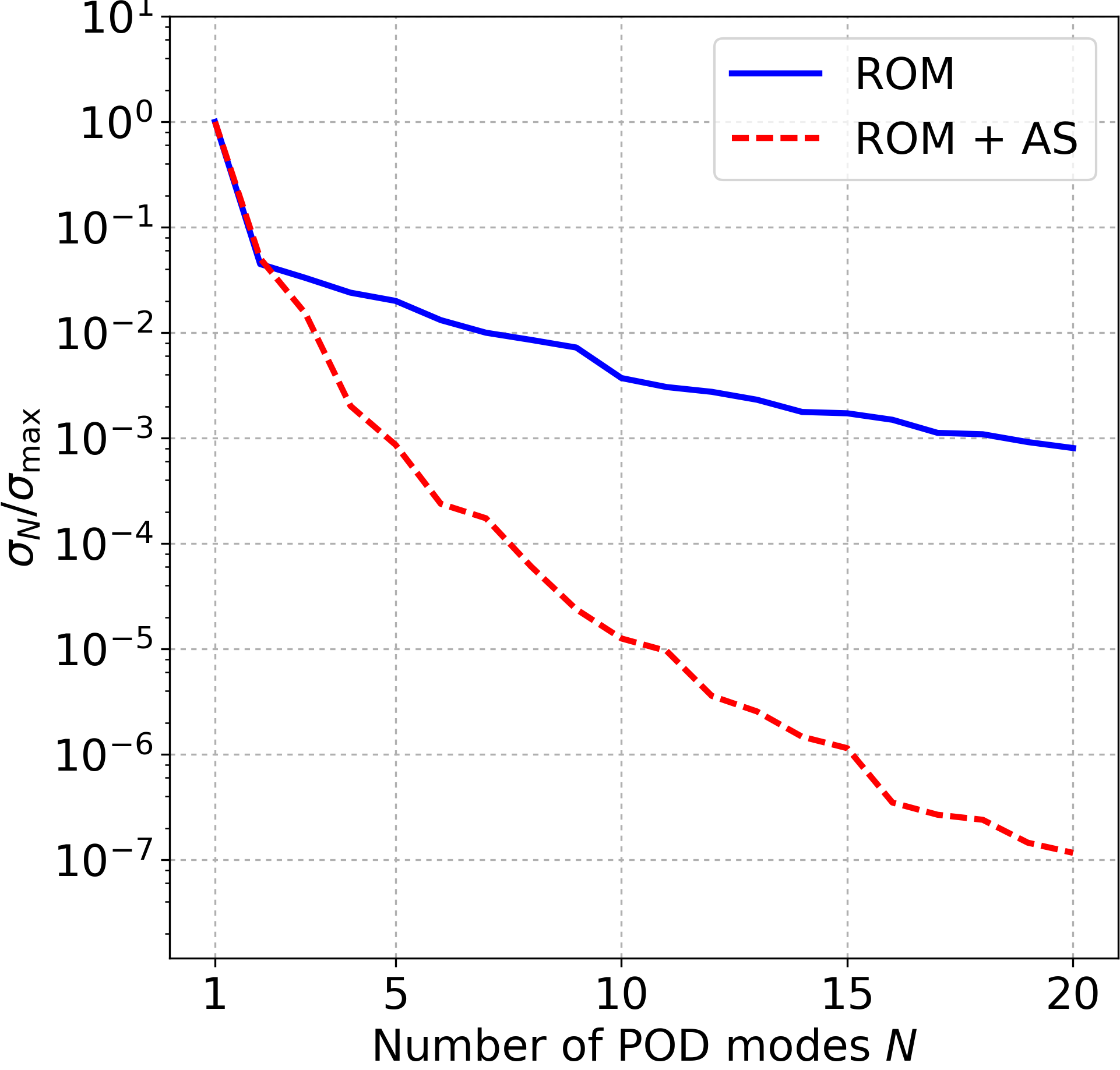}\label{subfig:press_sval}}
\caption{POD singular values as a function of the number $N$ of POD modes for%
\protect\subref{subfig:vel_sval} velocity, %
\protect\subref{subfig:sup_sval} supremizers, and %
\protect\subref{subfig:press_sval} pressure%
.}
\label{fig:singular_values}
\end{figure}

\begin{figure}[ht]
\centering
\subfloat[Velocity error.]{\includegraphics[width=0.48\textwidth]{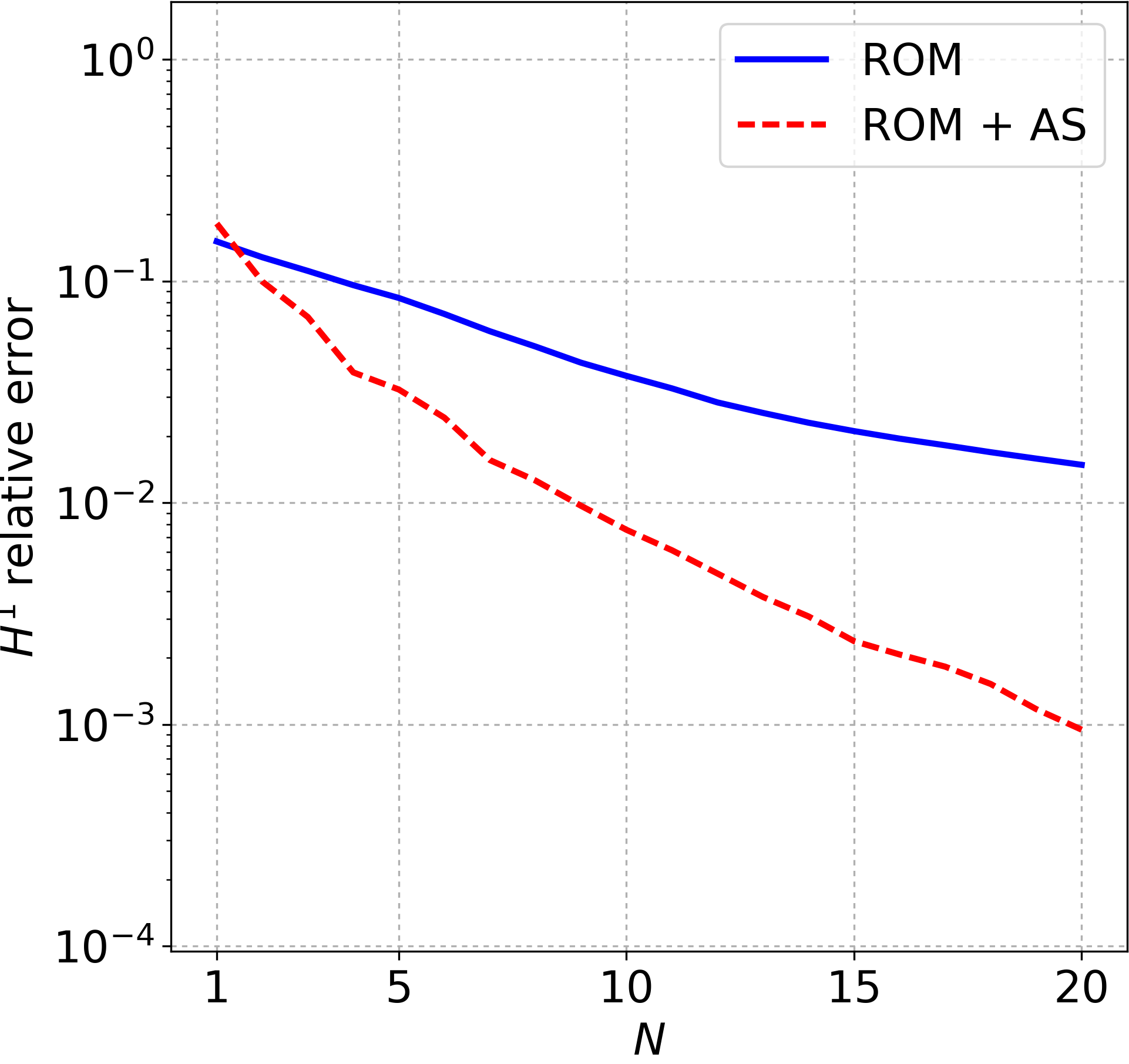}\label{subfig:vel_err}}\hfill
\subfloat[Pressure error.]{\includegraphics[width=0.48\textwidth]{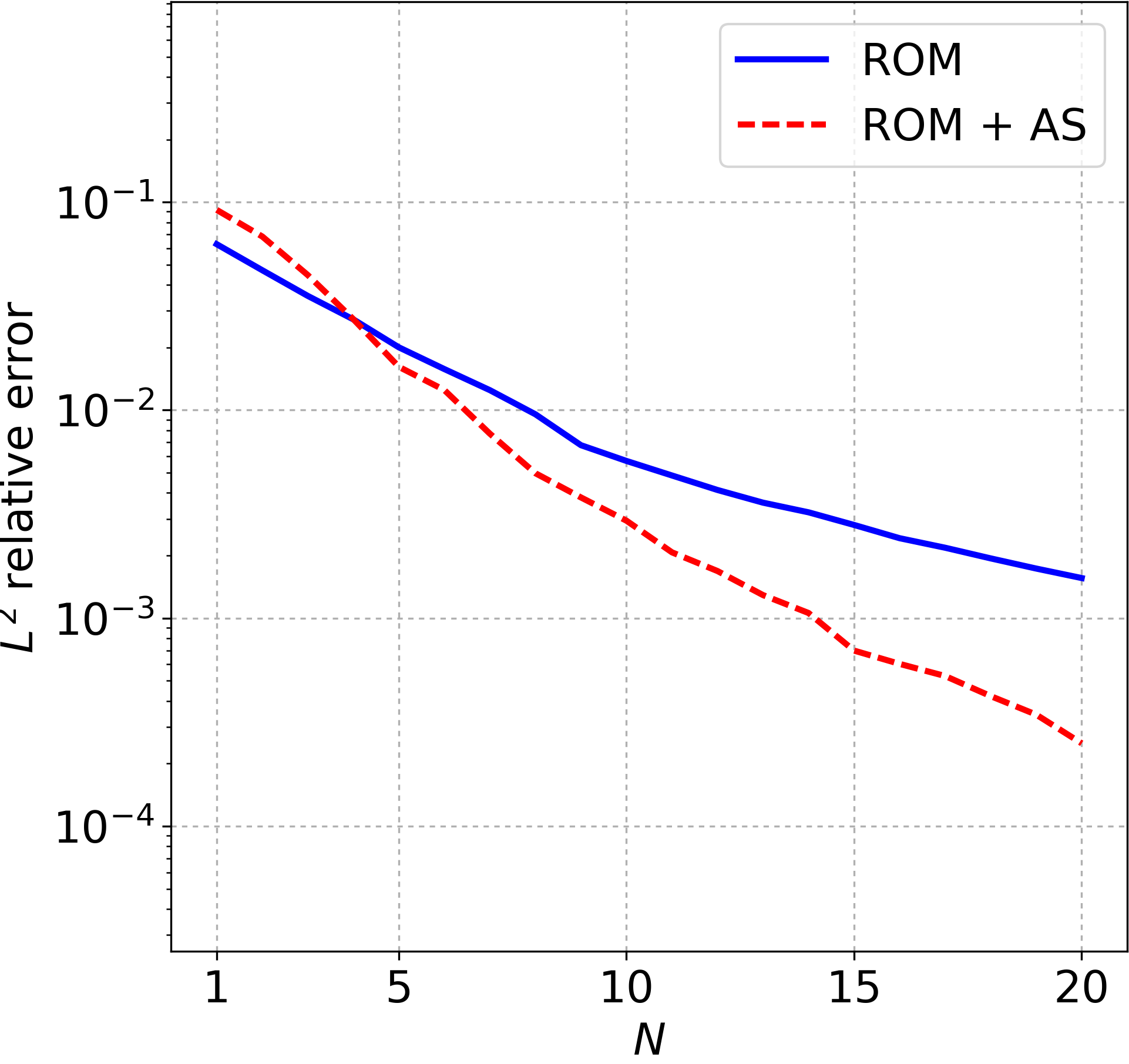}\label{subfig:press_err}}
\caption{Relative errors between the high fidelity solution and the reduced order one, as a function $N$, for \protect\subref{subfig:vel_err} velocity and \protect\subref{subfig:press_err} pressure component.}
\label{fig:error_as_rom}
\end{figure}

Once the active subspace $\mathbf{W}_1$ has been identified, we turn to the POD-Galerkin model order reduction defined in Section~\ref{sec:pod}.
The proposed combined methodology (denoted as ``ROM + AS'' in Figures \ref{fig:singular_values}-\ref{fig:error_as_rom}) will be compared to the standard POD-Galerkin approach on the full parameter space (denoted as ``ROM'' in Figures \ref{fig:singular_values}-\ref{fig:error_as_rom}) in order to highlight the effectiveness of our proposed method. 
The open source library RBniCS \cite{hesthaven2015certified,RBniCS} is employed to implement both methodologies.
During the training phase, we select a training set of size $N_{\text{train}}^{\text{POD}}$ and compute a POD of the resulting snapshots. In the case of the combined methodology we have chosen $N_{\text{train}}^{\text{POD}} = 100$, while in the standard approach we have chosen $N_{\text{train}}^{\text{POD}} = N_{\text{train}}^{\text{AS}}$. Corresponding POD singular values are show in Figure~\ref{fig:singular_values} for velocity, supremizers and pressure, as a function of the number $N$ of selected POD modes. The results show a slower 
decay for the standard approach when compared to the combined one,
meaning that the standard approach has to deal with a considerably
larger solution manifold. If the final goal is to evaluate the
quantity of interest, it is much more convenient to use the combined
method, which is able to provide a much smaller solution manifold by
neglecting the inactive (and so less interesting) directions. Indeed, Figure~\ref{fig:error_as_rom} shows that the combined methodology is able to
reach relative errors which are up to an order of magnitude smaller
when compared to the standard one, for both velocity
(Figure~\ref{subfig:vel_err}) and pressure
(Figure~\ref{subfig:press_err}) when $N = 20$. The errors are average relative errors on a testing set of cardinality $N_{\text{test}}^{\text{POD}} = 100$. A similar error analysis can be carried out for the quantity of interest, showing a trend similar to the one in Figure~\ref{subfig:press_err}.

\section{Conclusions and perspectives}
\label{sec:the_end}
In this chapter we have presented a combined parameter and computational model reduction by means of active subspaces and POD-Galerkin methods, and we applied the proposed combined method on a synthetic problem related to the estimation of the pressure drop across a stenosed artery bifurcation.
First, we reduced the high dimensional parameter space into a lower dimensional parameter subspace by means of the active subspaces property.
Our numerical test case, related to the deformation of two stenoses, shows an effective reduction of the dimensionality of the parameter space,
from 10 control points displacements to 2 active variables. Just the active parameter subspace is then employed for a further model reduction by means of a POD-Galerkin method. When comparing the performance (in terms of errors) of the resulting reduced order model with a standard one (without active subspaces preprocessing), the proposed approach shows better results up to an order of magnitude.
This is due to the fact that the standard approach has to account for several directions in the parameter space (the inactive subspaces) which only account for negligible variations in the pressure drop and thus could have been neglected. The proposed methodology could find further developments in more realistic cardiovascular problems, for what concerns both the geometry (e.g. patient's personalization) and the mathematical model (e.g. unsteadiness, compliance). Moreover, several enhancements on the combination of the two approaches could be foreseen; among the possible ones, we mention a more tight coupling between the training stages of active subspaces and POD-Galerkin sampling methods based on a greedy approach in order to avoid the solution of several high fidelity problems for two (possibly large) training sets.

\section*{Acknowledgements}
This work was partially supported by the INDAM-GNCS 2017 project ``Advanced numerical methods combined with computational reduction techniques
for parameterised PDEs and applications'', and by European Union Funding for Research and Innovation --- Horizon 2020 Program --- in the framework of European Research Council Executive Agency: H2020 ERC CoG 2015 AROMA-CFD project 681447 ``Advanced Reduced Order Methods with Applications in Computational Fluid Dynamics'' P.I. Gianluigi Rozza.


\end{document}